\documentclass[12pt]{amsart} 
\usepackage{amssymb,amsmath,amscd,mathrsfs} 
\usepackage{a4wide}
\usepackage{verbatim}
\usepackage[all]{xy} 
\usepackage{color}

\newtheorem{theorem}{Theorem}[section]

\theoremstyle{definition}

\theoremstyle{remark}

\newcommand{\cC}{\mathcal{C}}
\newcommand{\cI}{\mathcal{I}}
\newcommand{\cJ}{\mathcal{J}}

\newcommand{\cU}{\mathcal{U}}

\newcommand{\Fr}{\mathscr{F}}

\newcommand{\Hr}{\mathscr{H}}

\newcommand{\Dc}{\mathcal{D}}

\newcommand{\im}{\mathrm{Im}}

\newcommand{\C}{\mathbb{C}}

\newcommand{\Hp}{\mathbb{H}}
\newcommand{\R}{\mathbb{R}}

\newcommand{\Z}{\mathbb{Z}}

\newcommand{\rmi}{{\rm i}}
\def\cchi{\raisebox{.45 ex}{$\chi$}}
\def\build#1_#2^#3{\mathrel{\mathop{\kern 0pt#1}\limits_{#2}^{#3}}}
\begin{document} 
\title[On the instability of eigenvalues]{On the instability of
  eigenvalues}
\author{Sylvain Gol\'enia} 
\address{Institut de Math\'ematiques de Bordeaux, Universit\'e
Bordeaux $1$, $351,$ cours de la Libération
\\$33405$ Talence cedex, France}
\email{sylvain.golenia@u-bordeaux1.fr}
\begin{abstract} This is the proceeding of a 
talk given in Workshop on Differential Geometry and its applications 
at Alexandru Ioan Cuza University Ia\c si, Romania, September 2--4, 2009.
I explain how positive commutator estimates help
in the analysis of embedded eigenvalues in a geometrical setting.
Then, I will discuss the disappearance of eigenvalues in the
perturbation theory and its relation with the Fermi golden rule.
\end{abstract} 
\maketitle 
\section{Introduction}
Let $\Hp:=\{(x,y)\in \R^2, y>0\}$ be the Poincar\'e
half-plane and we endow it with the metric $g:=y^{-2}(dx^2+
dy^2)$. Consider the group $\Gamma:=PSL_2(\Z)$. It acts faithfully on $\Hp$ 
by homographies, from the left. The interior of a fundamental domain
of the quotient $\Hp \backslash \Gamma$ is given by $X:=\{(x,y)\in
\Hp, |x|<1, x^2+y^2>1\}$. Let $\Hr:=L^2(X, g)$ be the set of $L^2$
integrable function acting on $X$, with respect to the volume element
$dx\,dy/y^2$. Let $\cC^\infty_b(X)$ be the restriction to $X$ of the
smooth bounded functions acting on $\Hp$ which are $\C$-valued and invariant
under $\Gamma$. The (non-negative) Laplace operator is defined as the closure of
\begin{eqnarray*}
\Delta:=-y^2(\partial^2_x + \partial^2_y), \mbox{ on } \cC^\infty_b(X). 
\end{eqnarray*} 
It is a (unbounded) self-adjoint operator on $L^2(X)$. Using
Eisenstein series, for instance, one sees that its essential spectrum
is given by $[1/4, \infty)$ and that it has no singularly continuous
spectrum, with respect to the Lebesgue measure. It is well-known that
$\Delta$ has infinitely many eigenvalues accumulating at $+\infty$ and
that every eigenspace is of finite dimension. We refer to
\cite{F} for an introduction to the subject.

We consider the Schr\"odinger operator $H_\lambda:=\Delta+\lambda V$,
where $V$ is the multiplication by a bounded, real-valued function and
$k\in \R$. We focus on an eigenvalue $k>1/4$ of $\Delta$ and
assume that the following hypothesis of \emph{Fermi golden rule} holds
true. Namely, there is $c_0>0$ so that:
\begin{eqnarray}\label{e:FGR0}
\lim_{\varepsilon \rightarrow 0^+}
PV\overline{P}\,\im(H_0-k +\rmi\varepsilon)^{-1}
\overline{P}V P \geq c_0 P,
\end{eqnarray} 
in the form sense and where $P:=P_{k}$, the projection on the
eigenspace of $k$, and $\overline P:=1-P$. As $P$ is of finite
dimension, the limit can be taken in the weak or in the strong
sense. At least formally, $\overline{P}\,\im(H_0-k +\rmi\varepsilon)^{-1}
\overline{P}$ tends to the Dirac mass $\pi \delta_k(\overline P
H_0)$. Therefore, the potential $V$ couples   the eigenspace of $k$
and $\overline P H_0$ over $k$ in a non-trivial way. This is a key
assumption in the second-order perturbation theory of embedded
eigenvalues, e.g., \cite{RS2}, and all the art is to prove that it
implies there is $\lambda_0>0$ that $H_\lambda$ has no eigenvalue in a
neighborhood of $k$ for $\lambda\in (0, |\lambda_0|)$. 

In \cite{CdV}, one shows that generically the eigenvalues disappear
under the perturbation of a potential (or of the metric) on a compact set. In this note, we are interested about the optimal decay at infinity of the perturbation given by a potential. Using the general result obtained in \cite{CGH} and under a hypothesis of Fermi golden rule, one is only able to cover the assumption $VL^{3}=o(1)$, as $y\rightarrow +\infty$, where
$L$ denotes the operator of multiplication by $L:=(x,y)\mapsto 1+\ln(y)$.
We give the main result:

\begin{theorem}\label{t:main}
Let $k>1/4$ be an $L^2-$eigenvalue of $\Delta$. Suppose that
$VL=o(1)$, as $y\rightarrow +\infty$ and that the Fermi
golden rule \eqref{e:FGR0} holds true, then there is $\lambda_0>0$,
so that $H_\lambda$ has no eigenvalue in a neighborhood of $k$, for
all $\lambda\in (0, |\lambda_0|)$. Moreover, if $VL^{1+\varepsilon
}=o(1)$, as $y\rightarrow +\infty$ for some $\varepsilon>0$, then
$H_\lambda$ has no singularly continuous spectrum. 
\end{theorem} 

We believe that the hypothesis $VL=o(1)$ is optimal in the scale of $L$.
In our approach, we use the Mourre theory, see \cite{ABG, mou} and
establish a positive commutator estimate.

\section{Idea of the proof}
Standardly, for $y$ large enough
and up to some isometry $\cU$, see for instance \cite{FH, GMo1, GMo2} the
Laplace operator can be written as  
\begin{align}\label{e:H_0} 
\tilde \Delta= (-\partial^2_r + 1/4)\otimes
P_{0} + \tilde \Delta (1 \otimes P_0^\perp)
\end{align}
on $\cC^\infty_c\big((c,\infty), dr \big)\otimes \cC^\infty(S^1)$, for
some $c>0$ and where $P_0$ is the projection on constant functions and 
$P_0^\perp:=1-P_0$. The Friedrichs extension of the operator
$\tilde \Delta (1 \otimes P_0^\perp)$ has compact resolvent.

Then, as in \cite{GMo1, GMo2}, we construct a conjugate operator. 
One chooses $\Phi\in\cC^\infty_c(\R)$ with $\Phi(x)=x$ on $[-1,1]$, and sets
$\Phi_\Upsilon(x):=\Upsilon\Phi(x/\Upsilon)$, for $\Upsilon \geq
1$. Let $\tilde \cchi$ be a smooth cut-off
function being $1$ for $r$ big enough and $0$ for $r$ being close to $c$.
We define on  $\cC^\infty_c\big((c, \infty)\times S^1\big)$ a
micro-localized version of the generator of dilations:
\begin{equation}\label{e:SR}
S_{\Upsilon,0}:= \tilde\cchi \left(\big(\Phi_\Upsilon(-\rmi\partial_r) r + r
\Phi_\Upsilon(-\rmi\partial_r)\big)\otimes P_0 \,\right)\tilde\cchi.
\end{equation}
The operator $\Phi_\Upsilon(-\rmi \partial_r)$ is defined on the real line
by $\Fr^{-1}\Phi_\Upsilon(\cdot)\Fr$, where $\Fr$ is the unitary Fourier
transform.  We also denote its closure by $S_{\Upsilon,0}$ and it is self-adjoint. 
In \cite{FH} for instance, one does not use a
micro-localization and one is not able to deal with really singular
perturbation of the metric as in \cite{GMo1, GMo2}. 

Now, one obtains 
\begin{equation*}
[\partial_r^2, \tilde\cchi (\Phi_\Upsilon r +r
\Phi_\Upsilon ) \tilde\cchi]
=4\tilde\cchi\partial_r\Phi_\Upsilon\tilde\cchi +
\mbox{remainder.}
\end{equation*}
Using a cut-off function $\tilde \mu$ being $1$ on the
cusp and $0$ for $y\leq 2$, we set 
\begin{align}\label{e:SRvrai}
S_{\Upsilon}:= \cU^{-1}S_{\Upsilon,0}\,\cU\,\widetilde\mu
\end{align} 
This is self-adjoint in $L^2(X)$. Now by taking $\Upsilon$ big
enough, one can show, as in \cite{GMo1, GMo2} that 
given an interval $\cJ$ around $k$, there exist $\varepsilon_\Upsilon>
 0$ and a compact operator $K_\Upsilon$ such that the inequality
 \begin{equation}\label{e:mourre}
E_\cJ(\Delta) [\Delta,\rmi S_{\Upsilon}] E_\cJ(\Delta) \geq (4 \inf(\cJ)
 - \varepsilon_\Upsilon) E_\cJ(\Delta)+ E_\cJ(\Delta)K_\Upsilon
 E_\cJ(\Delta) \end{equation}  
holds in the sense of forms, and such that $\varepsilon_\Upsilon$ tends to $0$
as $\Upsilon$ goes to infinity. Here, $E_\cJ(\cdot)$ denotes the
spectral measure above the interval $\cJ$.  

Now, we apply $\overline{P}$ to the left and right of \eqref{e:mourre}. 
Easily one has $\overline{P}E_\cJ(\Delta)=
\overline{P}E_\cJ\big(\Delta \overline P\big)$. We get: 
 \begin{align*}
\overline{P} E_\cJ(\overline{P}\Delta) \big[\,\overline{P}\Delta,\rmi
\overline{P}S_{\Upsilon}\overline{P}\,\big]
E_\cJ(\overline{P}\Delta)\overline{P} \geq& \,\,(4 \inf(\cJ) 
 - \varepsilon_\Upsilon) \overline{P}E_\cJ(\Delta\overline{P})\overline{P}
 \\
 &\,\, + \overline{P}E_\cJ(\overline{P}\Delta)K_\Upsilon
 E_\cJ(\overline{P}\Delta) \overline{P} 
 \end{align*} 
One can show that $\overline{P}S_{\Upsilon}\overline{P}$ is
self-adjoint in $\overline{P} L^2(X)$ and that $\big[\,P\overline
\Delta, \overline{P}S_{\Upsilon}\overline{P}\,\big]$ extends to a
bounded operator.  

We now shrink the size of the interval $\cJ$. As
$\overline{P}\Delta$ has no eigenvalue in $\cJ$, then the operator
$\overline{P}E_\cJ(\overline{P}\Delta)K_\Upsilon
E_\cJ(\overline{P}\Delta) \overline{P}$ tends to $0$ in
norm. Therefore, by shrinking enough, one obtains a smaller interval
$\cJ$ containing $k$ and a constant $c>0$ so that 
\begin{align}\label{e:mourre3}
\overline{P} E_\cJ(\overline{P}\Delta) \big[\,\overline{P}\Delta,\rmi
\overline{P}S_{\Upsilon}\overline{P}\,\big]
E_\cJ(\overline{P}\Delta)\overline{P} \geq& \,\,c
\overline{P}E_\cJ(\Delta\overline{P})\overline{P} 
\end{align} 
holds true in the form sense on $\overline{P} L^2(X)$. At least
formally, the positivity on $\overline{P} L^2(X)$ of the commutator
$\big[\,H_\lambda ,\rmi \overline{P}S_{\Upsilon}\overline{P}\,\big]$,
up to some spectral measure and to some small $\lambda$, should be a
general fact and should not rely on the  
Fermi golden rule hypothesis.

We now try to extract some positivity on $P L^2(X)$. First, we set
\begin{eqnarray}\label{e:Rbar}
R_\varepsilon:=\big((H_0-k)^2+\varepsilon^2\big)^{-1/2},\,
\overline{R_\varepsilon}:=\overline{P} R_\varepsilon \mbox{ and }
F_\varepsilon := \overline{R_\varepsilon}^2. 
\end{eqnarray} 
Note that $\varepsilon R_\varepsilon^2=\im(H_0-k +\rmi\varepsilon)^{-1}$
and  that $R_\varepsilon$ commutes with $P$. Using \eqref{e:FGR0}, we get:  
\begin{eqnarray}\label{e:FGR}
(c_1/\varepsilon) P\geq
PV \overline{P}\, F_\varepsilon\, \overline{P} V P\geq (c_2/\varepsilon) P,  
\end{eqnarray} 
for $\varepsilon_0 > \varepsilon > 0$. 

We follow an idea of \cite{BFSS}, which was successfully used in
\cite{Go, Merkli} and set   
\begin{eqnarray*}
B_\varepsilon :=\im (\overline{R_\varepsilon}^2V P).
\end{eqnarray*} 
It is a finite rank operator. Observe now that we gain some positivity
as soon as $\lambda\neq 0$:
  \begin{eqnarray}\label{e:motiv}
P[H_\lambda, \rmi \lambda B_\varepsilon]P= \lambda^2 P V
F_\varepsilon V  P \geq (c_2\lambda^2 /\varepsilon) P.
\end{eqnarray}
It is therefore natural to modify the conjugate operator $S_\Upsilon$
to obtain some positivity on $PL^2(X)$. We set 
\begin{eqnarray}\label{e:Ahat}
\hat S_\Upsilon:=\overline P S_\Upsilon \overline P+\lambda \theta B_\varepsilon.
\end{eqnarray}  
It is self-adjoint on $\Dc(S_\Upsilon)$ and is diagonal with respect to
the decomposition $\overline PL^2(X)\oplus  P L^2(X)$. 

Here $\theta>0$ is a technical parameter. We choose   
 $\varepsilon$ and $\theta$, depending on $\lambda$, so that 
$\lambda=o(\varepsilon)$, $\varepsilon = o(\theta)$ and $\theta=o(1)$
as $\lambda$ tends to $0$. We summarize this into:
\begin{eqnarray}\label{e:size}
|\lambda|\,\, \ll\,\, \varepsilon \,\,\ll\,\, \theta\,\,\ll\,\, 1,
\mbox{ as } \lambda \mbox{  tends to } 0. 
\end{eqnarray} 

With respect to the  decomposition $\overline P E_\cJ(\Delta)\oplus  P
E_\cJ(\Delta)$, as $\lambda$ goes to $0$, we have  
\begin{align*}
E_\cJ(\Delta) \left[\lambda V, \rmi \overline P S_\Upsilon \overline P\,\right] E_\cJ(\Delta)&=
\left(\begin{array}{cc}
O(\lambda)& O(\lambda)\\
O(\lambda)& 0
\end{array}\right), 
\\
E_\cJ(\Delta)[\Delta, \rmi \lambda\theta B_\varepsilon ]E_\cJ(\Delta)&=
\left(\begin{array}{cc}
0& O(\lambda\theta\varepsilon^{-1/2})\\
O(\lambda\theta\varepsilon^{-1/2}) & 0
\end{array}\right),
\\
\mbox{ and }
E_\cJ(\Delta)[\lambda V, \rmi \lambda\theta B_\varepsilon]E_\cJ(\Delta)&=
\left(\begin{array}{cc}
O(\lambda^2\theta\varepsilon^{-3/2})& O(\lambda^2\theta\varepsilon^{-3/2})\\
O(\lambda^2\theta\varepsilon^{-3/2}) & \lambda^2\theta F_\varepsilon  
\end{array}\right).
\end{align*} 

Now comes the delicate point. Under the condition \eqref{e:size} and by
choosing $\cI$, slightly smaller than $\cJ$, we use the previous
estimates and a Schur Lemma to deduce:
\begin{equation}\label{e:mourrefinal}
E_\cI(H_\lambda) [H_\lambda,\rmi \hat S_{\Upsilon}] E_\cI(H_\lambda)
\geq \frac{c \lambda^2 \theta}{\varepsilon}E_\cI(H_\lambda), 
\end{equation} 
for some positive $c$ and as $\lambda$ tends to $0$.

We mention that only the decay of $VL$ is used to establish the
last estimate. In fact, one uses that $[V, \rmi \hat
S_\Upsilon](\Delta+1)^{-1}$ is a compact operator.  

Now it is a standard use of the Mourre theory to deduce Theorem
\ref{t:main} and refer to \cite{ABG}, see \cite{GMo1, GMo2} for some
similar application of the theory. For the absence of eigenvalue, one
relies on the fact that given an eigenfunction $f$ of $H_\lambda$
w.r.t.\ an eigenvalue $\kappa\in \cI$, one has: 
\begin{eqnarray}\label{e:virial}
\langle f, [H_\lambda,\rmi \hat S_{\Upsilon}] f\rangle = \langle f,
[H_\lambda -\kappa,\rmi \hat S_{\Upsilon}] f\rangle=0. 
\end{eqnarray}
Then, one applies $f$ on the right and on the left of
\eqref{e:mourrefinal} and infers that $f=0$ thanks to the fact that
the constant $c\lambda^2\theta$ is non-zero. 

In \cite{GMo1, GMo2}, we prove that the $C_0$-group $(e^{\rmi S_{\Upsilon}
  t})_{t\in\R}$ stabilizes the domain $\Dc(H_\lambda)=
\Dc(\Delta)$. By perturbation, we prove that this is also the case for
$(e^{\rmi \hat S_{\Upsilon}   t})_{t\in\R}$.  Thanks to this property,
we can expand the commutator of \eqref{e:virial} in a legal way. This
is known as the Virial theorem in the Mourre Theory, see 
\cite{ABG, mou}.

\bibliographystyle{plain} 
 
\end{document}